\documentclass[11pt,a4paper,twoside]{article}
\usepackage{amsmath,amsfonts,amssymb,amsthm}
\date{}
\newcommand{\Rn}{{\mathbb R}^n}
\newcommand{\D}{{\rm d}}

\begin{document}

\title{ A characterization of the mixed discriminant}
\renewcommand{\thefootnote}{\fnsymbol{footnote}}
\author{D. Florentin, V. D. Milman\footnote{Supported in part by Minerva Foundation, ISF grant 826/13 and by BSF grant 2012111.}, \hspace{2mm}R. Schneider} 

\maketitle

\renewcommand{\thefootnote}{{}}

\footnote{2010 {\em Mathematics Subject Classification.} Primary 52A39, 15A15

{\em Key words and phrases.} Mixed discriminant, mixed volume, characterization, $L_2$ addition, centered ellipsoid, polarization formula
}
\section{Introduction}

The striking analogy between mixed volumes of convex bodies and mixed discriminants of positive semidefinite matrices has repeatedly been observed. It was A. D. Aleksandrov \cite{Ale38} who, in his second proof of the Aleksandrov--Fenchel inequalities for the mixed volumes, first introduced the mixed discriminants of positive semidefinite quadratic forms and established some of their properties, including quadratic inequalities. This use of mixed discriminants in the theory of mixed volumes is described, for example, by Busemann \cite[Sect. 7]{Bus58} and Leichtwei{\ss} \cite{Lei93}.

In \cite{MS11}, the mixed volume of centrally symmetric convex bodies in $\Rn$ was characterized, up to a factor, as the only function of $n$ centrally symmetric convex bodies which is Minkowski additive and increasing (with respect to set inclusion) in each variable and which vanishes if two of its arguments are parallel segments. The strong analogy mentioned above leads one to expect a similar characterization of the mixed discriminant, but it appears that the arguments employed in \cite{MS11} cannot be transferred directly. The present note utilizes the $L_2$ addition of ellipsoids with center at the origin, represented by symmetric positive semidefinite matrices. The 
$L_2$ addition is a special case of the $L_p$ addition ($p\ge 1$) of convex bodies, which was introduced by Firey \cite{Fir62} and developed with great success by Lutwak, beginning with \cite{Lut93}. It turned out that the $L_2$ addition opens the way to employ geometric arguments similar to those used in \cite{MS11}, which allows us to obtain the desired characterization. It is formulated as Theorem 2 below.

\section{Statement of the result}

Let ${\mathcal M}^n$ denote the set of real symmetric positive semidefinite $n\times n$ matrices $(n\ge 1$). The {\em mixed discriminant} $D:({\cal M}^n)^n\to {\mathbb R}$ is the unique symmetric function for which
$$ {\det}(\lambda_1 A_1+\dots +\lambda_m A_m) = \sum_{i_1,\dots,i_n=1}^m \lambda_{i_1}\cdots\lambda_{i_n} D(A_{i_1},\dots,A_{i_n}) $$
for $m\in{\mathbb N}$, $A_1,\dots,A_m\in{\mathcal M}^n$ and $\lambda_1,\dots,\lambda_m\ge 0$. If $A^{(i)}$ denotes the $i$th column of the matrix $A$, then
$$ D(A_1,\dots,A_n)= \frac{1}{n!} \sum_{\sigma\in S(n)} \det(A^{(1)}_{\sigma(1)},\dots, A^{(n)}_{\sigma(n)}),$$
where $S(n)$ denotes the group of permutations of the numbers $1,\dots,n$. 

The mixed discriminant $D$ is a multilinear function, and it is nonnegative. We formulate a criterion for its positivity. To show the similarity to the conditions for the positivity of the mixed volume (\cite[Th. 5.1.7]{Sch93}), we denote by $P(A)$, for $A\in{\mathcal M}^n$, the subspace of ${\mathbb R}^n$ spanned by the eigenvectors of $A$ corresponding to positive eigenvalues.

\vspace{3mm}

\noindent{\bf Proposition 1.} {\em For $A_1,\dots,A_n\in{\mathcal M}^n$, the following assertions are equivalent:\\[1mm]
$\rm (a)$ $D(A_1,\dots,A_n)>0$,\\[1mm]
$\rm (b)$ there are linearly independent vectors $v_i\in P(A_i)$, $i=1,\dots,n$,\\[1mm]
$\rm (c)$ $\dim(P(A_{i_1})+\dots+P(A_{i_k})) \ge k$ for each choice of indices $1\le i_1<\dots<i_k\le n$ and for all $k\in\{1,\dots,n\}$,\\[1mm]
$\rm (d)$ If $K_i=P(A_i)\cap B^n$, where $B^n$ denotes the unit ball of $\Rn$, then
$$ V(K_1,\dots,K_n)>0$$
for the mixed volume $V$.}

\vspace{3mm}

The equivalence of (a) and (c) was proved by Panov \cite[Th. 1]{Pan87}, and the equivalence of (b) and (c) follows from \cite[Lemma 5.1.8]{Sch93}; see also Bapat \cite[Th. 9]{Bap89}. The equivalence with (d) is clear from \cite[Th. 5.1.7]{Sch93}.

We have formulated this criterion, first to stress once more the analogy between mixed volumes and mixed discriminants, and second because we shall need a rudimentary form of it in our characterization of the mixed discriminant.

\vspace{3mm}

\noindent{\bf Theorem 2.} {\em Let $F:({\mathcal M}^n)^n\to{\mathbb R}$ $(n\ge 1)$ be a nonnegative function which is additive in each variable and which is zero if two of its arguments are proportional matrices of rank one. Then there is a constant $a\ge 0$ with
$$ F(A_1,\dots,A_n)=aD(A_1,\dots,A_n)$$
for all $A_1,\dots,A_n\in{\mathcal M}^n$.}

\vspace{3mm}

We note that the assumption on the vanishing of $F$ is essential for the characterization. Without it, there are many functions 
$F:({\mathcal M}^n)^n\to{\mathbb R}$ which are multilinear and nonnegative. One example is given by
$$ F(A_1,\dots,A_n)= \int_{({\mathbb S}^{n-1})^n} \prod_{i=1}^n \langle u_i,A_iu_i\rangle \,\D\mu(u_1,\dots,u_n),$$
where ${\mathbb S}^{n-1}$ and $\langle\cdot,\cdot\rangle$ denote, respectively, the unit sphere and the scalar product  of $\Rn$ and where $\mu$ is a finite Borel measure on $({\mathbb S}^{n-1})^n$. Another example is 
$$ F(A_1,\dots,A_n)=\prod_{i=1}^n D(A_i,B_2,\dots,B_n),\quad A_1,\dots,A_n\in{\mathcal M}^n,$$
where $B_2,\dots,B_n\in{\mathcal M}^n$ are fixed. This function $F$ shares with the mixed discriminant the property of being symmetric (which we have not assumed in the theorem).

As in \cite{MS11}, we can also conclude that within the class of nonnegative functionals on ${\mathcal M}^n$ which are additive in each variable, the mixed discriminant is characterized by a minimality property.

\vspace{3mm}

\noindent{\bf Corollary 3.} {\em Let $F:({\mathcal M}^n)^n\to{\mathbb R}$ be a nonnegative function which is additive in each variable. If $F$ is bounded from above by a constant multiple of the mixed discriminant, then $F$ is itself a constant multiple of the mixed discriminant.}

\vspace{3mm}

Indeed, if $F(A_1,\dots,A_n)\le cD(A_1,\dots,A_n)$ and if two of the $A_i$'s are proportional matrices of rank one, then $D(A_1,\dots,A_n)=0$ and hence $F(A_1,\dots,A_n)=0$, so $F$ satisfies the assumptions of our theorem.

\vspace{3mm}

Let $F:({\mathcal M}^n)^n\to{\mathbb R}$ be additive in each variable. If $F$ is nonnegative, then $F$ is increasing in each variable with respect to the order $\prec$, where $A\prec B$ means that $B-A$ is positive semidefinite. Conversely, suppose that $F$ is increasing in at least one variable, say the first one. Denoting by $0_n$ the $n\times n$ zero matrix, we have $F(0_n,\dots) =F(0_n+0_n,\dots) = F(0_n,\dots)+F(0_n,\dots)$ (where the dots stand for the fixed remaining arguments) and hence $F(0_n,\dots)=0$. If $A\in{\mathcal M}^n$, then
$0_n\prec A$ and hence $F(A,\dots)\ge F(0_n,\dots)=0$. Thus $F$ is nonnegative.

\section{Centered ellipsoids and mixed discriminants}

For the proof of the theorem, we use the correspondence between positive semidefinite symmetric matrices and centered ellipsoids. Here an ellipsoid is called {\em centered} if its center is at the origin.

Recall that $B^n$ is the unit ball of ${\mathbb R}^n$. Let $T:{\mathbb R}^n\to{\mathbb R}^n$ be a linear mapping (possibly degenerate). Then $E= TB^n$ is an ellipsoid (not necessarily full-dimensional) which is centered. We denote the set of all centered ellipsoids in ${\mathbb R}^n$ by ${\mathcal E}^n$.

For a convex body $K\subset{\mathbb R}^n$, let $h(K,\cdot)$ denote its support function. For $u\in{\mathbb R}^n$  we have
\begin{eqnarray*}
h(E,u) &=& h(TB^n,u) =\max\{\langle x,u\rangle: x\in TB^n\} =\max\{\langle Tb,u\rangle: b\in B^n\}\\
&=& \max\{\langle b,T^*u\rangle: b\in B^n\}= h(B^n,T^*u)\\
&=& \|T^*u\|= \langle T^*u,T^*u\rangle^{1/2} =  \langle u, TT^*u\rangle^{1/2}.
\end{eqnarray*}

We use the standard orthonormal basis of ${\mathbb R}^n$ to identify an endomorphism of ${\mathbb R}^n$ with the matrix describing it with respect to this basis, and we interpret the vectors of $\Rn$ as columns, if necessary. Then $A= TT^*$ uniquely defines an element of ${\mathcal M}^n$. We have seen that
\begin{equation}\label{1} 
h(E,u)^2=\langle u,Au\rangle\quad\mbox{for }u\in\Rn.
\end{equation}
Conversely, each matrix $A\in{\mathcal M}^n$ determines a centered ellipsoid $E$ in this way. If (\ref{1}) holds, we write $A=A_E$ and we also say that $A$ and $E$ {\em correspond to each other}. 
Clearly,
$$ \dim E = {\rm rank}\, A_E.$$

The $L_2$ {\em sum} $E_1+_2 E_2$ of two centered ellipsoids $E_1,E_2$  is defined by
\begin{equation}\label{0}  
h(E_1+_2 E_2,\cdot)^2 = h(E_1,\cdot)^2+h(E_2,\cdot)^2,
\end{equation}
and the scalar multiple $\lambda\cdot_2 E$ by $h(\lambda\cdot_2 E,\cdot)^2= \lambda h(E,\cdot)^2$, thus $\lambda\cdot_2 E=\sqrt{\lambda} E$. Obviously, $E_1+_2 E_2$ is again a centered ellipsoid, and
$$ A_{E_1+_2 E_2}= A_{E_1}+A_{E_2},\quad A_{\lambda\cdot_2 E}= \lambda A_E.$$
It is important to notice that the $L_2$ sum does not depend on the dimension of the ambient space. Let $L$ be a subspace of $\Rn$ and let $E_1,E_2$ be centered ellipsoids contained in $L$. Then the $L_2$ sum $E_1+_2 E_2$ formed in $L$ is the same as the $L_2$ sum of $E_1$ and $E_2$ formed in ${\mathbb R}^n$. This follows from the fact that for a convex body $K\subset L$ and a vector $u\in{\mathbb R}^n$ we have $h(K,u)=h(K,u_L$), where $u=u_L+u_{L^\perp}$ with $u_L\in L$ and $u_{L^\perp}\in L^\perp$.

We add an observation which, though it is not needed to its full extent, is of general interest, as it allows to see mixed discriminants in even closer analogy to mixed volumes. 

For the volume ${\rm vol}$ we have
$$ {\rm vol}\, E = \kappa_n\sqrt{\det A_E},$$
where $\kappa_n$ denotes the volume of $B^n$. Hence, for $E_1,\dots,E_m\in{\mathcal E}^n$ and $\lambda_1,\dots,\lambda_m\ge 0$ we get
\begin{equation}\label{7} 
{\rm vol}(\lambda_1\cdot_2 E_1+_2\dots+_2 \lambda_m\cdot_2 E_m)^2= \kappa_n^2 \sum_{i_1,\dots,i_n=1}^m \lambda_{i_1}\cdots\lambda_{i_n}D(A_{i_1},\dots,A_{i_n}),
\end{equation}
where $A_i  = A_{E_i}$.

To emphasize the analogy, we recall that on the class of convex bodies, equipped with Minkowski linear combination and the size functional ${\rm vol}$, we have a polarization formula leading to the mixed volume (see, e.g., \cite{Sch93}, (5.1.16)). Equation (\ref{7}) shows that on the class of centered ellipsoids, equipped with $L_2$ linear combination and the size functional ${\rm vol}^2/\kappa_n^2$, a polarization formula yields the mixed discriminant of positive semidefinite symmetric matrices.

\section{Proof of the characterization theorem}

Let $F:({\mathcal M}^n)^n\to{\mathbb R}$ be a nonnegative function which is additive in each variable and which is zero if two of its arguments are proportional matrices of rank one. We define the mapping $G:({\mathcal E}^n)^n\to{\mathbb R}$ by 
\begin{equation}\label{5} 
G(E_1,\dots,E_n)= F(A_{E_1},\dots,A_{E_n}).
\end{equation} 
Then $G$ is additive in each variable with respect to the $L_2$ sum $+_2$, and it is zero if two of its arguments are parallel centered segments. We assert that $G$ is increasing under set inclusion in each of its arguments. In fact, let $E, E', E_2,\dots,E_n\in{\mathcal E}^n$ be such that $E\subset E'$ and let $A,A',A_2,\dots,A_n\in{\mathcal M}^n$ be the corresponding matrices. For all $u\in{\mathbb R}^n$ we have $h(E,u)\le h(E',u)$, hence $\langle u,Au\rangle\le\langle u,A'u\rangle$. Thus, the matrix $A'-A$ is positive semidefinite. Since the mapping $F$ is nonnegative, this gives $F(A'-A,A_2,\dots,A_n)\ge 0$. Since $F$ is additive in its first argument, we obtain $F(A,A_2,\dots,A_n)\le F(A',A_2,\dots,A_n)$ and thus $G(E,E_2,\dots,E_n)\le G(E',E_2,\dots,E_n)$. Here the first argument can be replaced by any other argument.

Next, we assert that
\begin{equation}\label{4}
G(\lambda\cdot_2 E,E_2,\dots,E_n)=\lambda G(E,E_2,\dots,E_n)
\end{equation}
for $\lambda\ge 0$ and $E,E_1,\dots,E_n\in{\mathcal E}^n$, and similarly in the other arguments. In fact, for $k\in {\mathbb N}$ it follows from (\ref{0}) that
$$ \underbrace{E+_2\dots+_2 E}_{k \mbox{ times}} =\sqrt{k}E,$$
hence the additivity of $G$ gives $G(\sqrt{k}E,\dots)= kG(E,\dots)$. This leads to $G(\sqrt{q} E,\dots)=qG(E,\dots)$ for rational $q>0$. If $\lambda>0$ is any real number, we choose rational numbers $p,q>0$ with $p<\lambda<q$, and from $\sqrt{p}E\subset\sqrt{\lambda} E \subset \sqrt{q}E$ and the monotonicity of $G$ we obtain $pG(E,\dots)=G(\sqrt{p} E,\dots)\le G(\sqrt{\lambda} E,\dots)\le G(\sqrt{q} E, \dots)= qG(E,\dots)$. With $p,q\to\lambda$ this yields $G(\sqrt{\lambda} E,\dots)=\lambda G(E,\dots)$ and thus (\ref{4}) for $\lambda>0$. Part of the preceding argument also shows that $G(0\cdot_2 E,\dots)=0$.

The following arguments are modelled after those of \cite{MS11}, with some essential modifications. Here and below we write ${\mathbb R}S=\{\lambda x: \lambda\in{\mathbb R},\,x\in S\}$.

\vspace{3mm}

\noindent{\bf Lemma 4.} {\em Let $S,T_1,T_2\subset{\mathbb R}^n$ be nondegenerate centered segments satisfying
\begin{equation}\label{2} 
T_1+{\mathbb R}S=T_2+{\mathbb R}S.
\end{equation}
Then 
$$ G(T_1,S,E_1,\dots,E_{n-2})= G(T_2,S,E_1,\dots,E_{n-2})$$
for all $E_1,\dots,E_{n-2}\in{\mathcal E}^n$.

Here the role of the first two arguments can be played by any two other arguments, in any order.}

\vspace{3mm}

\noindent{\em Proof.} By (\ref{2}), the segments $S,T_1,T_2$ are contained in a two-dimensional subspace $L$. Let $0<\varepsilon<1$ be given. Let $\pm v\in L$ be the unit vectors orthogonal to $S$. We may assume that none of the segments $T_1,T_2$ is parallel to $S$, since otherwise the assertion is trivial. Then, by (\ref{2}), $h(T_1,v)=h(T_2,v)>0$ and, therefore, $\sqrt{1-\varepsilon}\,h(T_1,v)< h(T_2,v)$. There is a neighborhood $U$ of $v$ with $U=-U$ such that
\begin{equation}\label{3} 
\sqrt{1-\varepsilon}\, h(T_1,u)< h(T_2,u)\quad\mbox{for all } u\in L\cap {\mathbb S}^{n-1}\cap U.
\end{equation}
There is a number $\alpha>0$ with $h(S,u)\ge \alpha$ for all $u\in L\cap {\mathbb S}^{n-1}\setminus U$. Therefore, we can choose a number $\lambda_\varepsilon>0$ with
$$ (1-\varepsilon)h(T_1,u)^2\le h(T_2,u)^2+\lambda_\varepsilon h(S,u)^2\quad\mbox{for all } u\in L\cap {\mathbb S}^{n-1}\setminus U.$$
By (\ref{3}), this holds also for $u\in L\cap {\mathbb S}^{n-1}\cap U$. Since the convex sets $T_1$ and $T_2+_2 \lambda_\varepsilon\cdot_2 S $ are contained in $L$, the latter inequality for the support functions shows that
$$ (1-\varepsilon)\cdot_2 T_1\subset T_2+_2 \lambda_\varepsilon\cdot_2 S$$
holds in $L$ and therefore also if $+_2$ refers to the ambient space $\Rn$. Since $G$ is increasing under set inclusion, we obtain
$$ G((1-\varepsilon)\cdot_2 T_1,S,E_1,\dots,E_{n-2}) \le G(T_2+_2 \lambda_\varepsilon\cdot_2 S,S,E_1,\dots,E_{n-2})$$
and thus, since $G$ is $\cdot_2$-homogeneous and $+_2$-additive in each argument,
\begin{eqnarray*} 
(1-\varepsilon) G(T_1,S,E_1,\dots,E_{n-2}))&\le& G(T_2,S,E_1,\dots,E_{n-2})+ \lambda_\varepsilon G(S,S,E_1,\dots,E_{n-2})\\
&=& G(T_2,S,E_1,\dots,E_{n-2}),
\end{eqnarray*}
where we have used that $G$ vanishes if two of its arguments are parallel segments. Now we can let $\varepsilon$ tend to zero and obtain
$$ G(T_1,S,E_1,\dots,E_{n-2}) \le G(T_2,S,E_1,\dots,E_{n-2}).$$
Since $T_1$ and $T_2$ can be interchanged, the assertion follows. \qed

\vspace{3mm}

\noindent{\bf Lemma 5.} {\em $G$ is symmetric on centered segments.}

\vspace{3mm}

\noindent{\em Proof.} Let $S,T\subset {\mathbb R}^n$ be centered segments of positive length. If they are parallel, then $G(S,T,\dots)=0=G(T,S,\dots)$. Suppose they are not parallel. Let $D$ be a diagonal of the parallelogram $S+T$. Since $D+{\mathbb R}S =T +{\mathbb R}S$, Lemma 4 gives $G(D,S,\dots)=G(T,S,\dots)$. Since $D+{\mathbb R}T =S +{\mathbb R}T$, Lemma 4 also gives $G(D,T,\dots)=G(S,T,\dots)$. Since $S+{\mathbb R}D=T+{\mathbb R}D$, Lemma 4 (with arguments in the different order) gives $G(D,S,\dots)=G(D,T,\dots)$. Altogether we obtain $G(S,T,\dots)=G(T,S,\dots)$. Similarly, the symmetry of $G$ in any other pair of arguments is obtained. \qed

Now we are in a position to prove the main theorem.

For $n=1$, $F$ is a real function on the nonnegative real numbers which is additive and hence a solution of Cauchy's functional equation. Since it is nonnegative, it must be of the form $F(x)=ax$, $x\ge 0$, with a nonnegative constant $a$. This proves the case $n=1$ of the theorem.

Now we assume that $n\ge 2$ and that the assertion of the theorem has been proved in lower dimensions. Let $S$ be a nondegenerate centered segment and let $S^\perp$ be the linear subspace of $\Rn$ orthogonal to $S$. We first assume that $S$ is parallel to the vector $e_n$ of the standard basis $(e_1,\dots,e_n)$ of $\Rn$, so that $(e_1,\dots,e_{n-1})$ is an orthonormal basis of $S^\perp$.

For $E_2,\dots,E_n\in{\mathcal E}^n$ with $E_i\subset S^\perp$ we define
$$ g(E_2,\dots,E_n)= G(S,E_2,\dots,E_n),$$
where $G$ is defined by (\ref{5}). Then $g$ is $+_2$-additive in each variable and nonnegative, and it vanishes if two of its arguments are proportional centered segments. For a centered ellipsoid $E\subset S^\perp$ we denote by $A'_E$ its matrix with respect to the basis $(e_1,\dots,e_{n-1})$. The mixed discriminant in $S^\perp$ is denoted by $D'$. Then it follows from the inductional hypothesis, applied to the function $(A'_{E_2},\dots, A'_{E_n})\mapsto g(E_2,\dots,E_n)$, that
$$ g(E_2,\dots,E_n)= c(S) D'(A'_{E_2},\dots,A'_{E_n}),$$
with a constant $c(S)$ depending on $S$.

With respect to the basis $(e_1,\dots,e_n)$, the matrix $A_S$ has entry $\ell(S)^2$, the square of the length of $S$, at position $(n,n)$ and zero at all other positions. For an ellipsoid $E\in{\mathcal E}^n$ with $E\subset S^\perp$, the matrix $A_E$ has zero at positions $(n,j)$ and $(i,n)$, and the remaining submatrix is given by $A'_E$. It follows that $\det(A_S+A_E)=\ell(S)^2\det(A'_E)$. From this we get by mixing (i.e., replacing $A_S$ by $\lambda_1A_S$ and $A_E$ by $\lambda_2A_{E_2}+\dots+\lambda_nA_{E_n}$ with $\lambda_i\ge 0$, expanding, and comparing the coefficients of $\lambda_1\cdots\lambda_n$) that
$$ nD(A_S,A_{E_2},\dots,A_{E_n})=\ell(S)^2D'(A'_{E_2},\dots,A'_{E_n}), $$
hence
\begin{equation}\label{6}
G(S,E_2,\dots,E_n)=a(S)D(A_S,A_{E_2},\dots,A_{E_n})
\end{equation}
with $a(S)= nc(S)/\ell(S)^2$. Since $G$ is positively $\cdot_2$-homogeneous and $D$ is positively homogeneous, we have $a(\lambda\cdot_2 S)=a(S)$ for $\lambda>0$.

For any centered ellipsoids $E_1,\dots,E_n\in{\mathcal E}^n$, the corresponding matrices $A_1,\dots,A_n$ with respect to the standard basis and the matrices $\bar A_1,\dots,\bar A_n$ with respect to a rotated image of this basis are related by $\bar A_i=B^{-1}A_iB$ with a suitable orthogonal matrix $B$. From $\det(B^{-1}AB)=\det A$ it follows that 
$$D(B^{-1}A_1B,\dots,B^{-1}A_nB)=D(A_1,\dots,A_n).$$
Therefore, the relation (\ref{6}) holds also without the special assumption on the direction of $S$ made above.

Now let $S,T_2,\dots,T_n$ be nondegenerate centered segments. Let $S_2$ be the image of $T_2$ under orthogonal projection to $S^\perp$. By Lemma 4, 
$$ G(S,T_2,T_3,\dots,T_n)= G(S,S_2,T_3,\dots,T_n).$$
This holds also if $S_2$ is degenerate, sin´ce then $T_2$ is parallel to $S$ and both sides are zero. Treating the remaining arguments similarly, we arrive at
\begin{eqnarray*} 
G(S,T_2,\dots,T_n)&=&G(S,S_2,\dots,S_n) = a(S)D(A_S,A_{S_2},\dots,A_{S_n})\\
&=&a(S)D(A_S,A_{T_2},\dots,A_{T_n}),
\end{eqnarray*}
where in the last step we have used that the function $D$ has the same properties as $F$, so that $S_i$ and $T_i$ can be interchanged, $i=2,\dots,n$. 

From Lemma 5 we get
$$ G(S,T_2,T_3,\dots,T_n)=G(T_2,S,T_3,\dots,T_n)= a(T_2)D(A_{T_2},A_S,A_{T_3},\dots,A_{T_n}),$$
from which we can conclude that $a(S)=a$ does not depend on $S$. 

Since each ellipsoid $E\in{\mathcal E}^n$ is a finite $+_2$-sum of centered segments, multilinearity can be used to show that
$$ G(E_1,\dots,E_n)= aD(A_{E_1},\dots,A_{E_n})$$
holds for $E_1,\dots,E_n\in{\mathcal E}^n$. This completes the proof. \qed

{\sc D. Florentin and V. D. Milman, School of Mathematical Sciences, Tel Aviv University, Ramat Aviv, Tel Aviv 69978, Israel

R. Schneider, Mathematisches Institut, Albert-Ludwigs-Universit\"at, Eckerstr. 1, D--79104, Freiburg i. Br., Germany}

\end{document}